\documentclass{article}
\usepackage{amsmath,amssymb,amsthm}

\newtheorem{theorem}{Theorem}
\newtheorem{proposition}{Proposition}
\newtheorem{definition}{Definition}
\newtheorem{corollary}{Corollary}

\numberwithin{theorem}{section}
\numberwithin{definition}{section}
\numberwithin{corollary}{section}
\numberwithin{equation}{section}
\numberwithin{proposition}{section}

\begin{document}

\title{Weak Order Complexes}
	\author{Sergei~Ovchinnikov\\
	Mathematics Department\\
	San Francisco State University\\
	San Francisco, CA 94132\\
	sergei@sfsu.edu}
	\date{\empty}
\maketitle

\begin{abstract}
\noindent
The paper presents geometric models of the set {\bf WO} of weak orders on a finite set $X$. In particular, {\bf WO} is modeled as a set of vertices of a cubical subdivision of a permutahedron. This approach is an alternative to the usual representation of {\bf WO} by means of a weak order polytope.
\end{abstract}

\section{Introduction}

Let $\mathcal{B}$ be a family of binary relations on a finite set $X$. This set can be endowed with various structures which are important in applications. One particular way to represent $\mathcal{B}$ is to embed it into a cube $\{0,1\}^N$ of sufficiently large dimension ($N=|X|^2$ would always work) by using characteristic functions of relations in $\mathcal{B}$, and consider a convex hull of the set of corresponding points. Then $\mathcal{B}$ is treated as a polytope with rich combinatorial and geometric structures. There are many studies of \emph{linear order polytopes}, \emph{weak order polytopes}, \emph{approval--voting polytopes}, and \emph{partial order polytopes}, and their applications. (See, for instance,~\cite{jD02,sF03,sFpF03} and references there.)

In this paper we study the set {\bf WO} of all weak orders on $X$ from a different point of view. Namely, we model the Hasse diagram of {\bf WO} as a $1$--skeleton of a cubical subdivision of a permutahedron. Our motivation has its roots in media theory~\cite{jF97,jFsO02,sOaD00} where it is shown that the graph of a medium is a partial cube~\cite{sOaD00}.

Section 2 presents some basic facts about weak orders and the Hasse diagram of {\bf WO}. In Section~3 we describe various geometric models of {\bf WO}. They are combinatorially equivalent under the usual connection between zonotopes, polar zonotopes, and hyperplane arrangements. Finally, in Section~4, we give an application of our approach to media theory by constructing a weak order medium.

\section{The Hasse diagram {\bf WO}}

In the paper, $X$ denotes a finite set with $n>1$ elements. A binary relation $W$ on $X$ is a \emph{weak order} if it is transitive and strongly complete. Antisymmetric weak orders are \emph{linear orders}. The set of all weak orders (resp. linear orders) on $X$ will be denoted {\bf WO} (resp. {\bf LO}).

For a weak order $W$, the \emph{indifference} relation $I=W\cap W^{-1}$ is an equivalence relation on $X$. Equivalence classes of $I$ are called \emph{indifference} classes of $W$. These classes are linearly ordered by the relation $W\slash I$. We will use the notation $W=(X_1,\ldots,X_k)$ where $X_i$'s are indifference classes of $W$ and $(x,y)\in W$ if and only if $x\in X_i,\;y\in X_j$ for some $1\leq i\leq j\leq k$. Thus our notation reflects the linear order induced on indifference classes by $W$.

We distinguish weak orders on $X$ by the number of their respective indifference classes: if $W=(X_1\ldots,X_k)$, we say that $W$ is a \emph{weak $k$--order}. The set of all weak $k$--orders will be denoted $\textbf{WO}(k)$. In particular, weak $n$--orders are linear orders and there is only one weak $1$--order on $X$, namely, $W=(X)=X\times X$, which we will call a \emph{trivial} weak order. Weak $2$--orders play an important role in our constructions. They are in the form $W=(A,X\setminus A)$ where $A$ is a nonempty proper subset of $X$. Clearly, there are $2^n-2$ distinct weak $2$--orders on a set of cardinality $n$.

The set {\bf WO} is a partially ordered set with respect to the set inclusion relation $\subseteq$. We denote the Hasse diagram of this set by the same symbol {\bf WO}. The following figure shows, as an example, {\bf WO} for a $3$--element set $X=\{a,b,c\}$.

\begin{center}
\begin{picture}(220,140)(-110,-70)

\multiput(-100,20)(40,0){6}{\circle*{8}}\multiput(-100,-20)(40,0){6}{\circle*{8}}\multiput(-100,-20)(40,0){5}{\line(1,1){40}}
\multiput(-100,-20)(40,0){6}{\line(0,1){40}}
\put(-100,20){\line(5,-1){200}}

\put(0,60){\circle*{8}}

\put(-100,20){\line(5,2){100}}
\put(-60,20){\line(3,2){60}}
\put(-20,20){\line(1,2){20}}
\put(20,20){\line(-1,2){20}}
\put(60,20){\line(-3,2){60}}
\put(100,20){\line(-5,2){100}}

\put(-95,-60){Fig. 1. The Hasse diagram of {\bf WO} for $n=3$.}

\end{picture}
\end{center}
Here the maximal element corresponds to the trivial weak order, the six vertices in the layer below correspond to weak $2$--orders, and the vertices in the lowest layer correspond to the linear orders on $X$.

We find it more intuitive to represent the Hasse diagram {\bf WO} by a directed graph as shown in the following figure. (Similar diagrams were introduced in~\cite[ch.2]{jKjS72} and~\cite{kB73}).

\begin{center}
\begin{picture}(220,220)(-110,-110)
\put(30,45){\circle*{5}}\put(-60,0){\circle*{5}}\put(0,-90){\circle*{5}}\put(0,90){\circle*{5}}\put(60,0){\circle*{5}}
\put(90,45){\circle*{5}}\put(-90,45){\circle*{5}}\put(90,-45){\circle*{5}}\put(-90,-45){\circle*{5}}
\put(30,-45){\circle*{5}}\put(-30,45){\circle*{5}}\put(-30,-45){\circle*{5}}\put(0,0){\circle*{5}}

\put(30,45){\vector(-2,-3){29}}\put(-30,-45){\vector(2,3){28}}\put(-30,45){\vector(2,-3){28}}
\put(30,-45){\vector(-2,3){29}}\put(60,0){\vector(-1,0){58}}\put(-60,0){\vector(1,0){57}}

\put(90,45){\vector(-1,0){58}}\put(-90,45){\vector(1,0){57}}\put(-90,-45){\vector(1,0){57}}\put(90,-45){\vector(-1,0){58}}
\put(0,90){\vector(2,-3){28}}\put(90,45){\vector(-2,-3){29}}\put(-90,-45){\vector(2,3){28}}\put(-90,45){\vector(2,-3){28}}
\put(90,-45){\vector(-2,3){29}}\put(0,-90){\vector(2,3){28}}\put(0,-90){\vector(-2,3){28}}\put(0,90){\vector(-2,-3){29}}

\put(4,88){\shortstack{$a$\\$b$\\$c$}}\put(94,36){\shortstack{$b$\\$a$\\$c$}}\put(-100,36){\shortstack{$a$\\$c$\\$b$}}
\put(94,-55){\shortstack{$b$\\$c$\\$a$}}\put(-100,-55){\shortstack{$c$\\$a$\\$b$}}\put(4,-109){\shortstack{$c$\\$b$\\$a$}}

\put(7,3){$abc$}
\put(28,48){\shortstack{$ab$\\$c$}}\put(-40,48){\shortstack{$a$\\$bc$}}\put(68,-4){\shortstack{$b$\\$ca$}}
\put(-82,-4){\shortstack{$ca$\\$b$}}\put(28,-64){\shortstack{$bc$\\$a$}}\put(-40,-64){\shortstack{$c$\\$ab$}}
\end{picture}
\end{center}

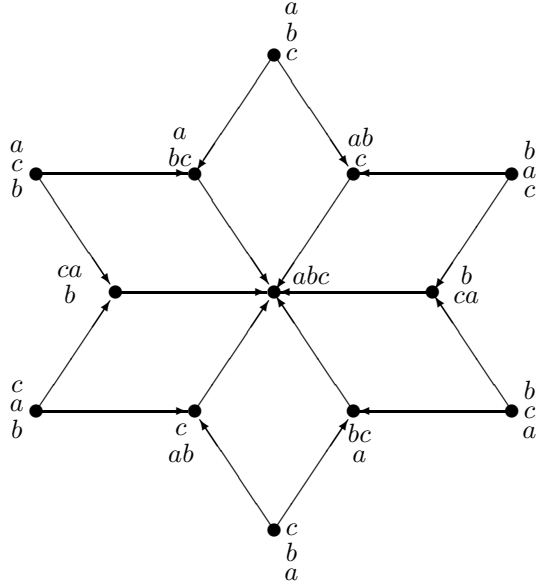
\centerline{Fig. 2. Another form of the Hasse diagram of {\bf WO} for $n=3$.}

\medskip\noindent
Here the arrows indicate the partial order on $\textbf{WO}$ and, for instance, the weak order $(\{ab\},\{c\})$ is represented as \shortstack{$ab$\\$c$}.

In the rest of this section we establish some properties of {\bf WO}. The following proposition corrects the statement of Problem 19 on p.115 in~\cite{bM79}.

\begin{proposition} \label{SuborderStructure}
A weak order $W'$ contains a weak order $W=(X_1,\ldots,X_k)$ if and only if
\begin{equation*}
W'=\left(\bigcup_{j=1}^{i_1} X_j,\bigcup_{j=i_1+1}^{i_2} X_j,\ldots,\bigcup_{j=i_m}^k X_j\right)
\end{equation*}
for some sequence of indices $1\leq i_1<i_2\cdots<i_m\leq k$.
\end{proposition}

\begin{proof}
Let $W\subset W'$. Then the indifference classes of $W$ form a subpartition of the partition of $X$ defined by the indifference classes of $W'$. Thus each indifference class of $W'$ is a union of some indifference classes of $W$. Since $W\subset W'$, we can write $W'=(\cup_1^{i_1} X_j,\cup_{i_1+1}^{i_2} X_j,\ldots,\cup_{i_m}^k X_j)$ for some sequence of indeces $1\leq i_1<\cdots<i_m\leq k$.

\end{proof}
One can say~\cite[ch.2]{bM79} that $W\subset W'$ if and only if the indifference classes of $W'$ are ``enlargements of the adjacent indifference classes'' of $W$.

\begin{corollary} \label{Cover}
A weak order $W'$ covers a weak order $W=(X_1,\ldots,X_k)$ in the Hasse diagram {\bf WO} if and only if $W'=(X_1,\ldots,X_i\cup X_{i+1},\ldots,X_k)$ for some $1\leq i<k$.
\end{corollary}

\begin{proposition} \label{UniqueRepresentation}
A weak order admits a unique representation as an intersection of weak $2$--orders, i.e., for any $W\in\text{\emph{\bf WO}}$ there is a uniquely defined set $J\subseteq\text{\emph{\bf WO}}(2)$ such that
\begin{equation} \label{IntersectionRepresentation}
W = \bigcap_{U\in J} U.
\end{equation}
\end{proposition}

\begin{proof}
Clearly, the trivial weak order has a unique representation in the form (\ref{IntersectionRepresentation}) with $J=\emptyset$.

Let $W=(X_1,\ldots,X_k)$ with $k>1$ and let $J_W$ be the set of all weak $2$--orders containing $W$. By Proposition~\ref{SuborderStructure}, each weak order in $J_W$ is in the form
\begin{equation*} \label{form}
W_i = (\cup_{1}^i X_j,\cup_{i+1}^k X_j),\quad 1\leq i<k.
\end{equation*}
Let $(x,y)\in\bigcap_{i=1}^{k-1} W_i$. Suppose $(x,y)\notin W$. Then $x\in X_p$ and $y\in X_q$ for some $p>q$. It follows that $(x,y)\notin W_q$, a contradiction. This proves (\ref{IntersectionRepresentation}) with $J=J_W$.

Let $W$ be a weak order in the form (\ref{IntersectionRepresentation}). Clearly, $J\subseteq J_W$. Suppose that $W_s=(\cup_1^s X_1,\cup_{s+1}^k X_j)\notin J$ for some $s$. Let $x\in X_{s+1}$ and $y\in X_s$. Then $(x,y)\in W_i$ for any $i\not=s$, but $(x,y)\notin W$, a contradiction. Hence, $J=J_W$ which proves uniqueness of representation (\ref{IntersectionRepresentation}). 

\end{proof}

Let $J_W$, as in the above proof, be the set of all weak $2$--orders containing $W$, and let $\mathcal{J}=\{J_W\}_{_{W\in\text{\bf WO}}}$ be the family of all such subsets of $\text{\bf WO}(2)$. The set $\mathcal{J}$ is a poset with respect to the inclusion relation.

The following theorem is an immediate consequence of Proposition~\ref{UniqueRepresentation}.

\begin{theorem} \label{DualIsomorphism}
The correspondence $W\mapsto J_W$ is a dual isomorphism of posets {\bf WO} and $\mathcal{J}$.
\end{theorem}

Clearly, the trivial weak order on $X$ corresponds to the empty subset of $\text{\bf WO}(2)$ and the set {\bf LO} of all linear orders on $X$ is in one--to--one correspondence with maximal elements in $\mathcal{J}$. The Hasse diagram {\bf WO} is dually isomorphic to the Hasse diagram of $\mathcal{J}$.

\begin{theorem} \label{SimplicialComplex}
The set $\mathcal{J}$ is a combinatorial simplicial complex, i.e., $J\in\mathcal{J}$ implies $J'\in\mathcal{J}$ for all $J'\subseteq J$.
\end{theorem}

\begin{proof}
Let $J'\subseteq J=J_W$ for some $W\in\text{\bf WO}$, i.e., $W=\bigcap_{_{U\in J_W}} U$. Consider $W'=\bigcap_{_{U\in J'}} U$. Clearly, $W'$ is transitive. It is complete, since $W\subseteq W'$. By Proposition~\ref{UniqueRepresentation}, $J'=J_{W'}\in\mathcal{J}$.

\end{proof}

It follows that $\mathcal{J}$ is a complete graded meet--semilattice. Therefore the Hasse diagram {\bf WO} is a complete join--semilattice with respect to the join operation $W\lor W' = \overline{W\cup W'}$, the transitive closure of $W\cup W'$.

\section{Geometric models of {\bf WO}}

A weak order polytope $\text{\bf P}_{WO}^n$ is defined as the convex hull in $\mathbb{R}^{n(n-1)}$ of the characteristic vectors of all weak orders on $X$ (see, for instance,~\cite{sFpF03}). Here we suggest different geometric models for {\bf WO}. For basic definitions in the area of polytopes and complexes, the reader is referred to Ziegler's book~\cite{gZ95}.

\begin{definition}
A \emph{cube} is a polytope combinatorially equivalent to $[0,1]^m$. A \emph{cubical complex} is a polytopal complex $\mathcal{C}$ such that every $P\in\mathcal{C}$ is a cube. The \emph{graph} $G(\mathcal{C})$ of a cubical complex $\mathcal{C}$ is the $1$--skeleton of $\mathcal{C}$.
\end{definition}

Thus the vertices and the edges of $G(\mathcal{C})$ are the vertices and the edges of cubes in $\mathcal{C}$, and $G(\mathcal{C})$ is a simple undirected graph.

Let $d=2^n-2$, where $n=|X|$, be the number of elements in $\text{\bf WO}(2)$. We represent each $W\in\text{\bf WO}$ by a characteristic function $\chi(J_W)$ of the set $J_W$. These characteristic functions are vertices of the cube $[0,1]^d$. Let $L\in\text{\bf LO}$ be a linear order on $X$. Then $J_L$ is a maximal element in $\mathcal{J}$ and, by Theorem~\ref{SimplicialComplex}, the convex hull of $\{\chi(J_W)\}_{_{W\supseteq L}}$ is a subcube $C_L$ of $[0,1]^d$. The dimension of $C_L$ is $n-1$. The collection of all cubes $C_L$ with $L\in\text{\bf LO}$ and all their subcubes form a cubical complex $\mathcal{C}(\text{\bf WO})$ which is a subcomplex of $[0,1]^d$. Clearly, $\mathcal{C}(\text{\bf WO})$ is a pure complex of dimension $n-1$ and the graph of this complex is isomorphic to the graph (that we denote by the same symbol, {\bf WO}) of the Hasse diagram of {\bf WO}.  

The above construction yields an isometric embedding of the graph {\bf WO} into the graph of $[0,1]^d$. Thus the graph {\bf WO} is a partial cube. We will use this fact in the last section.

The dimension $\dim \mathcal{C}(\text{\bf WO})=n-1$ is much smaller than the dimension $d=2^n-2$ of the space $\mathbb{R}^d$ in which $\mathcal{C}(\text{\bf WO})$ was realized. Simple examples indicate that $\mathcal{C}(\text{\bf WO})$ can be realized in a space of a much smaller dimension. 
\begin{center}
\begin{picture}(200,200)(-100,-100)

\thicklines
\put(-70,0){\vector(1,0){140}}\put(75,-2){$x$}
\put(0,-70){\vector(0,1){140}}\put(-2,75){$z$}
\put(-40,-40){\vector(1,1){80}}\put(43,43){$y$}

\thinlines
\put(0,0){\circle*{4}}
\put(50,0){\circle*{4}}\put(-50,0){\circle*{4}}
\put(0,50){\circle*{4}}\put(0,-50){\circle*{4}}
\put(26,26){\circle*{4}}\put(-25,-25){\circle*{4}}
\put(75,25){\circle*{4}}\put(-75,-25){\circle*{4}}
\put(26,76){\circle*{4}}\put(-25,-75){\circle*{4}}
\put(50,-50){\circle*{4}}\put(-50,50){\circle*{4}}

\put(50,0){\line(1,1){25}}\put(26,26){\line(1,0){50}}
\put(26,26){\line(0,1){51}}\put(0,50){\line(1,1){25}}
\put(-50,0){\line(0,1){50}}\put(-50,50){\line(1,0){50}}
\put(-50,0){\line(-1,-1){25}}\put(-25,-25){\line(-1,0){50}}
\put(-25,-25){\line(0,-1){50}}\put(0,-50){\line(-1,-1){25}}
\put(50,0){\line(0,-1){50}}\put(0,-50){\line(1,0){50}}

\multiput(0,0)(2,2){13}{\line(1,0){50}}
\multiput(0,0)(2,2){13}{\line(0,1){50}}
\multiput(0,0)(-2,0){25}{\line(0,1){50}}
\multiput(0,0)(-2,0){25}{\line(-1,-1){25}}
\multiput(0,0)(0,-2){25}{\line(-1,-1){25}}
\multiput(0,0)(0,-2){25}{\line(1,0){50}}

\put(-50,-100){Fig. 3. ``Monkey Saddle''.}
\end{picture}
\end{center}

For instance, for $n=3$ we have a realization of $\mathcal{C}(\text{\bf WO})$ in $\mathbb{R}^3$ as shown in Figure~3. (This is a `flat' analog of the popular smooth surface $z=x^3-3xy^2$.) One can compare this picture with the picture shown in Figure~2.

It turns out that there is a cubical complex, which is combinatorially equivalent to $\mathcal{C}(\text{\bf WO})$, and such that its underlying set is a polytope in $\mathbb{R}^{n-1}$.

We begin with a simple example. Let $X=\{1,2,3\}$ and let $\Pi_2$ be the $2$--dimensional permutahedron. Consider a subdivision of $\Pi_2$ shown in Figure~4.

\begin{center}
\begin{picture}(220,220)(-110,-110)
\put(0,90){\circle*{5}} \put(-7,95){$123$}\put(90,45){\circle*{5}} \put(88,50){$132$}
\put(90,-45){\circle*{5}} \put(88,-56){$312$}\put(0,-90){\circle*{5}} \put(-7,-102){$321$}
\put(-90,45){\circle*{5}} \put(-104,-56){$231$}\put(-90,-45){\circle*{5}} \put(-104,50){$213$}

\put(0,0){\circle*{5}}\put(90,0){\circle*{5}}\put(-90,0){\circle*{5}}\put(45,-67.5){\circle*{5}}
\put(-45,67.5){\circle*{5}}\put(45,67.5){\circle*{5}}\put(-45,-67.5){\circle*{5}}

\put(0,90){\line(2,-1){90}}\put(90,45){\line(0,-1){90}}\put(90,-45){\line(-2,-1){90}}
\put(0,-90){\line(-2,1){90}}\put(-90,-45){\line(0,1){90}}\put(-90,45){\line(2,1){90}}

\put(-45,-67.5){\line(2,3){90}}\put(45,-67.5){\line(-2,3){90}}\put(-90,0){\line(1,0){180}}

\put(-100,-120){Fig. 4. A cubical complex associated with $\Pi_2$.}
\end{picture}
\end{center}

\bigskip\noindent
Clearly, this subdivision defines a cubical complex which is combinatorially isomorphic to the cubical complex shown in Figure~3. (Compare it also with the diagram in Figure~2.)

In general, let $\Pi_{n-1}$ be a permutahedron of dimension $n-1$, where $n=|X|$. According to~\cite[p.18]{gZ95}, ``$k$--faces (of $\Pi_{n-1}$) correspond to ordered partitions of (the set $X$) into $n-k$ nonempty parts'' (see also~\cite{mBBm70}, p.54). In other words, each face of $\Pi_{n-1}$ represents a weak order on $X$. Linear orders on $X$ are represented by the vertices of $\Pi_{n-1}$ and the trivial weak order on $X$ is represented by $\Pi_{n-1}$ itself. Weak $2$--orders are in one--to--one correspondence with the facets of $\Pi_{n-1}$. Let $L$ be a vertex of $\Pi_{n-1}$. Consider the set of barycenters of all faces of $\Pi_{n-1}$ containing $L$. A direct computation shows that the convex hull $C_L$ of these points is a (combinatorial) cube. This is actually true for any simple zonotope ($\Pi_{n-1}$ is a simple zonotope). The following argument belongs to G\"{u}nter Ziegler~\cite{gZ03}.

Let $Z$ be a simple zonotope. By Corollary~7.18 in~\cite{gZ95}, $C_L$ is the intersection of the vertex cone of $L$ (which is a simplicial cone) with the dual facet cone of the dual of $Z$ (which is again a simplicial cone). This intersection is an $(n-1)$--dimensional (combinatorial) cube.

Cubes in the form $C_L$ form a subdivision of $\Pi_{n-1}$ and, together with their subcubes, form a cubical complex isomorphic to $\mathcal{C}(\text{\bf WO})$.

Another geometric model for the set {\bf WO} of all weak orders on $X$ can be obtained using the polar polytope $\Pi_{n-1}^\Delta$. Let $L(\Pi_{n-1})$ be the face lattice of the permutahedron $\Pi_{n-1}$. The joint--semilattice {\bf WO} is isomorphic to the joint--semilattice $L(\Pi_{n-1})\setminus\{\emptyset\}$ (Figure~1). By duality, the Hasse diagram {\bf WO} is dually isomorphic to the meet--semilattice $L(\Pi_{n-1}^\Delta)\setminus\{\Pi_{n-1}^\Delta\}$ of all proper faces of $\Pi_{n-1}^\Delta$. Under this isomorphism, the linear orders on $X$ are in one--to--one correspondence with facets of $\Pi_{n-1}^\Delta$, the weak $2$--orders on $X$ are in one--to--one correspondence with vertices of $\Pi_{n-1}^\Delta$, and the trivial weak order on $X$ corresponds to the empty face of $\Pi_{n-1}^\Delta$. Note that $\Pi_{n-1}^\Delta$ is a simplicial polytope. The set of its proper faces is a simplicial complex which is a geometric realization of the combinatorial simplicial complex $\mathcal{J}$ (cf. Theorem~\ref{SimplicialComplex}).

Other geometric and combinatorial models of {\bf WO} can be constructed by using the usual connections between zonotopes, hyperplane arrangements, and oriented matroids~\cite{gZ95}. One particular model utilizes the following well known facts about weak orders on $X$.

Let $f$ be a real--valued function on $X$ and, as before, let $n=|X|$. Then $W_f$ defined by
\begin{equation*}
(x,y)\in W_f\quad\Leftrightarrow\quad f(x)\leq f(y),
\end{equation*}
for all $x,y\in X$, is a weak order. On the other hand, for a given weak order $W$ there exists a function $f$ such that $W=W_f$. Two functions $f$ and $g$ are said to be equivalent if $W_f=W_g$. Clearly, equivalent functions form a cone $C_W$ in $\mathbb{R}^n$ and the union of these cones is $\mathbb{R}^n$. Thus there is a natural one--to-one correspondence between the set {\bf WO} and the family $\{C_W\}_{_{W\in\text{\bf WO}}}$. The cones in the form $C_W$ arise from a hyperplane arrangement $\mathcal{H}$ defined by the hyperplanes $H_{ij}=\{x\in\mathbb{R}^n : x_i=x_j\}$. The arrangement $\mathcal{H}$ is the hyperplane arrangement associated with the zonotope $\Pi_{n-1}$. Following the standard steps~\cite{gZ95}, one can also construct an oriented matroid representing {\bf WO}.

Geometric objects introduced in this section, the cubical complex $\mathcal{C}(\text{\bf WO})$, the simplicial complex $\mathcal{J}$ of proper faces of the polar zonotope $\Pi_{n-1}^\Delta$, and the hyperplane arrangement $\mathcal{H}$, all share the combinatorial structure of the Hasse diagram {\bf WO}.

\section{Weak order media}

In this section we construct a medium having {\bf WO} as a set of states.

A \emph{medium} is a pair $(\mathcal{V,T})$, where $\mathcal{V}$ is a set, whose elements are called \emph{states}, and $\mathcal{T}$ is a set of functions mapping $\mathcal{V}$ into itself satisfying certain axioms; elements of $\mathcal{T}$ are called \emph{tokens}. Models based on media are natural tools in studies of preference evolution~\cite{jF96,jFjD97,jFmRbG97} and panel data~\cite{mRjFbG99}.

The reader is referred to~\cite{jF97,jFsO02,sOaD00} for formal definitions and main results in the area of media theory. It suffices to consider one particular example of a medium to understand the developments in this section.

Let $Z$ be a finite set. The \emph{distance} $d(A,B)$ between two sets $A,B\subseteq Z$ is defined by $d(A,B)=|A\Delta B|$. A family $\mathcal{F}$ of subsets of $Z$ is said to be \emph{well graded}~\cite{jDjF97} if, for any $A,B\in\mathcal{F}$ there is a sequence $A_0=A,A_1,\ldots,A_k=B$ of sets in $\mathcal{F}$ such that $d(A,B)=k$ and $d(A_i,A_{i+1})=1$ for all $0\leq i<k$. We assume that $\cup\mathcal{F}=\emptyset$ and $\cap\mathcal{F}=Z$. A combinatorial simplicial complex is an example of a well graded set.

A well graded family $\mathcal{F}$ is representable as a medium $(\mathcal{F,T})$, where $\mathcal{T}$ contains, for all $x\in Z$, the two transformations $\tau_x,\,\tilde{\tau}_x$ of $\mathcal{F}$ into $\mathcal{F}$ defined by, for $S\in\mathcal{F}$,
\begin{align}
&\tau_x:S\mapsto S\tau_x = \begin{cases} \label{token}
		S\cup\{x\} &\text{if $S\cup\{x\}\in\mathcal{F}$,} \\
		S &\text{otherwise.}
	\end{cases} \\
&\tilde{\tau}_x:S\mapsto S\tilde{\tau}_x = \begin{cases} \label{reverse}
		S\setminus\{x\} &\text{if $S\setminus\{x\}\in\mathcal{F}$,} \\
		S &\text{otherwise.}
	\end{cases}
\end{align}

Actually, in some precise sense, any medium is isomorphic to a medium of well graded sets~\cite{sOaD00}.

A token $\tau_x$ is \emph{effective} on $S$ if $S\tau_x=S'\not=S$. In this case $S'\tilde{\tau}_x=S$ and we say that $\tau_x$ and $\tilde{\tau}_x$ are mutual \emph{reverses}.

By Theorems~\ref{DualIsomorphism} and~\ref{SimplicialComplex}, the set {\bf WO} can be identified with elements of the well graded family $\mathcal{J}$ of subsets of $\text{\bf WO}(2)$. Therefore {\bf WO} is representable as a medium. Tokens in this medium are defined by weak $2$--orders according to~(\ref{token}) and~(\ref{reverse}).

In what follows we describe effective actions of tokens in $\mathcal{T}$ in terms of weak orders.

Let $V$ be a weak $2$--order. One can treat $V$ as a partition of the set $X$ of ``alternatives'' into the subset $G$ of ``good'' alternatives and the subset $B$ of ``bad'' alternatives in the representation $V=(B,G)$. We consider effective actions of $\tau_V$ and $\tilde{\tau}_V$ separately.

(i)~The case of $\tau_V$. Let $W=(X_1,\ldots,X_k)$ be a weak order on which $\tau_V$ is effective, i.e., $W\tau_V\not=W$. Then $V\notin J_W$ and $J_W\cup\{V\}\in\mathcal{J}$. Thus, by Proposition~\ref{UniqueRepresentation}, $W\cap V$ is a weak order and $W$ covers $W\cap V$ in {\bf WO}. The indifference classes of $W\cap V$ are in the forms $X_i\cap B$ and $X_i\cap G$. Since $W\not=W\cap V$, there is $p$ such that $X_p\cap B$ and $X_p\cap G$ form a partition of $X_p$. By Corollary~\ref{Cover},
\begin{equation*}
W\tau_V = (X_1,\ldots,X_p\cap B,X_p\cap G,\ldots,X_k),
\end{equation*}
where $X_i\subset B$ for $i<p$, $X_i\subset G$ for $i>p$. In other words, the action of $\tau_V$ partitions the indifference class $X_p$ of $W$ into subsets of ``bad'' and ``good'' alternatives according to $V$. The remaining indifference classes of $W$ consist entirely of either ``bad'' or ``good'' alternatives.

(ii)~The case of $\tilde{\tau}_V$. Again, let $W=(X_1,\ldots,X_k)$ be a weak order on which $\tilde{\tau}_V$ is effective. Then $V\in  J_W$ which implies $W\subset V$. By Proposition~\ref{SuborderStructure}, there is $1\leq p<k$ such that $V=(\cup_1^p X_i,\cup_{p+1}^k X_i)$. By Proposition~\ref{UniqueRepresentation}, the weak order $W\tilde{\tau}_V$ is the intersection of weak $2$--orders in the form $W_q=(\cup_1^q X_i,\cup_{q+1}^k X_i)$ with $q\not=p$. The indifference classes of $W\tilde{\tau}_V$ are intersections of the indifference classes of $W_q$'s with $q\not=p$. It follows that
\begin{equation*}
W\tilde{\tau}_V = (X_1,\ldots,X_p\cup X_{p+1},\ldots,X_k).
\end{equation*}
Since $B=\cup_1^p X_i$ and $G=\cup_{p+1}^k X_i$, the action of $\tilde{\tau}_V$ joins the indifference classes $X_p$ and $X_{p+1}$ consisting of maximal (with respect to $W$) ``bad'' alternatives and minimal ``good'' alternatives, respectively.

The actions of tokens $\tau_V$ and $\tilde{\tau}_V$ are illustrated in the following diagrams.
\begin{center}
\begin{picture}(300,100)(-150,-50)

\put(-145,35){$W$}
\put(-110,35){[}
\put(-108,37){\line(1,0){40}}
\put(-69,35){]}
\put(-95,28){$X_1$}
\multiput(-50,35)(5,0){3}{$\cdot$}
\put(-20,35){[}
\put(-18,37){\line(1,0){40}}
\put(21,35){]}
\put(-3,28){$X_p$}
\multiput(40,35)(5,0){3}{$\cdot$}
\put(70,35){[}
\put(72,37){\line(1,0){40}}
\put(111,35){]}
\put(87,28){$X_k$}

\put(-145,0){$V$}
\put(-110,0){[}
\put(-108,2){\line(1,0){110}}
\put(1,0){]}
\put(-60,-7){$B$}
\put(5,0){[}
\put(6,2){\line(1,0){105}}
\put(110,0){]}
\put(55,-7){$G$}

\put(-145,-40){$W\tau_V$}
\put(-110,-40){[}
\put(-108,-38){\line(1,0){40}}
\put(-69,-40){]}
\put(-95,-47){$X_1$}
\multiput(-50,-40)(5,0){3}{$\cdot$}
\put(-20,-40){[}
\put(-18,-38){\line(1,0){20}}
\put(1,-40){]}
\put(5,-40){[}
\put(6,-38){\line(1,0){18}}
\put(23,-40){]}
\put(-15,-47){$X_p'$}
\put(8,-47){$X_p''$}
\multiput(43,-40)(5,0){3}{$\cdot$}
\put(70,-40){[}
\put(72,-37){\line(1,0){40}}
\put(111,-40){]}
\put(87,-47){$X_k$}

\end{picture}
\end{center}

\centerline{$X_p=X_p'\cup X_p''$}

\medskip
\centerline{Fig.~5. An effective action of $\tau_V$.}

\medskip
\begin{center}
\begin{picture}(300,100)(-150,-50)

\put(-145,35){$W\tilde{\tau}_V$}
\put(-110,35){[}
\put(-108,37){\line(1,0){40}}
\put(-69,35){]}
\put(-95,28){$X_1$}
\multiput(-50,35)(5,0){3}{$\cdot$}
\put(-20,35){[}
\put(-18,37){\line(1,0){40}}
\put(20,35){]}
\put(-18,28){$X_p\cup X_{p+1}$}
\multiput(40,35)(5,0){3}{$\cdot$}
\put(70,35){[}
\put(72,37){\line(1,0){40}}
\put(110,35){]}
\put(87,28){$X_k$}

\put(-145,0){$V$}
\put(-110,0){[}
\put(-108,2){\line(1,0){110}}
\put(1,0){]}
\put(-60,-7){$B$}
\put(5,0){[}
\put(6,2){\line(1,0){105}}
\put(110,0){]}
\put(55,-7){$G$}

\put(-145,-40){$W$}
\put(-110,-40){[}
\put(-108,-38){\line(1,0){40}}
\put(-69,-40){]}
\put(-95,-47){$X_1$}
\multiput(-50,-40)(5,0){3}{$\cdot$}
\put(-20,-40){[}
\put(-18,-38){\line(1,0){20}}
\put(1,-40){]}
\put(5,-40){[}
\put(6,-38){\line(1,0){18}}
\put(23,-40){]}
\put(-15,-47){$X_p$}
\put(8,-47){$X_{p+1}$}
\multiput(43,-40)(5,0){3}{$\cdot$}
\put(70,-40){[}
\put(72,-37){\line(1,0){40}}
\put(111,-40){]}
\put(87,-47){$X_k$}

\end{picture}
\end{center}

\medskip
\centerline{Fig.~6. An effective action of $\tilde{\tau}_V$.}

\section{Concluding remarks}

\begin{enumerate}
\item Melvin Janowitz has noted to the author that a lattice theoretical study of the joint--semilatice {\bf WO} was presented in~\cite{mJ84}. In particular, it is shown there (Proposition~{\bf F1}) that the intervals above atoms of {\bf WO} are isomorphic to the lattice of all subsets of an $(n-1)$--element set. These intervals are exactly the cubes $C_L$ (Section~3).
\item Using standard numerical representations for semiorders, interval orders, and biorders, one can construct geometric models based on hyperplane arrangements for the families of these relations like it was done at the end of Section~3 for the family of weak orders. These models can be used, for instance, to prove, in a rather transparent way, the wellgradedness property of these families which was established in~\cite{jDjF97}.
\end{enumerate}

\section*{Acknowledgments}
An earlier version of this paper was presented at the \emph{Ordinal and Symbolic Data Analysis Conference} (OSDA 2003), Irvine, August 2003. The author wishes to thank participants of this conference for their valuable comments.
Special thanks go to Jean--Paul Doignon and Jean--Claude Falmagne for their careful reading of the original draft of the paper and Melvin Janowitz for the reference to his paper~\cite{mJ84}.

\end{document}